\documentclass[12pt]{article}
\usepackage{amsmath,amssymb}
\usepackage{pstricks}
\usepackage{pst-node} 
\usepackage{latexsym}
\usepackage{color,graphicx}
\usepackage{wrapfig}
\usepackage{caption3} 
\DeclareCaptionOption{parskip}[]{} 
\usepackage[small]{caption}
\usepackage{setspace}
\newtheorem{Lemma}{Lemma}
\newtheorem{Proposition}[Lemma]{Proposition}

\newtheorem{Conjecture}[Lemma]{Conjecture}

\newtheorem{Corollary}[Lemma]{Corollary}

\newcommand{\forn}{\mbox{${\cal F}_0$}}
\newcommand{\foro}{\mbox{${\cal F}_1$}}

\newcommand{\eps}{\varepsilon}

\newcommand{\bZ}{{\bf Z}}

\newcommand{\bX}{{\bf X}}

\newcommand{\sfrac}[2]{{\textstyle\frac{#1}{#2}}}

\newcommand{\bY}{{\bf Y}}
\newcommand{\bx}{{\bf x}}

\newcommand{\bc}{{\bf c}}

\newcommand{\MC}{mul\-ti\-pli\-ca\-tive co\-a\-le\-scent}

\newcommand{\cvd}{l^2_{\mbox{{\footnotesize $\searrow$}}}}

\def\endpf{$\Box$}
\def\endrem{$\diamond$}

\newcommand{\len}{\rm len}

\begin{document}
\title{A playful note on spanning and surplus edges}
\date{}
\author{Vlada Limic\thanks{
        UMR 8628, Département de Math\'ematiques,
        CNRS et Universit\'e Paris-Sud XI,
        91405 Orsay, France}\\
        e-mail: vlada.limic@math.u-psud.fr}

\maketitle
\begin{abstract}
Consider a (not necessarily near-critical) random graph running in continuous time.
A recent breadth-first-walk construction is extended in order to account for the surplus edge data in addition to the spanning edge data. 
Two different graph representations of the {\MC}, with different advantages and drawbacks,  are discussed in detail.
A canonical multi-graph of Bhamidi, Budhiraja and Wang (2014) naturally emerges. 
The presented framework should facilitate understanding of scaling limits with surplus edges for near-critical random graphs in the domain of attraction of general (not necessarily standard) eternal \MC.
\end{abstract}



\section{Introduction}
If $n\in \mathbb{N}$, write $[1,n]$ for $\{1,\ldots, n\}$.
A continuous time variation of the Erd\"os-R\'enyi \cite{erdren} random graph $G(n,p)$, where $p\in[0,1]$, is naturally constructed as follows: fix $n$ vertices $[1,n]$ and let each edge (out of ${n \choose 2}$) appear at rate $1$, independently of each other. 
Here and in the rest of the paper {\em connected} means connected by a path of edges in the usual graph theory sense. 
If the minimal path is in fact an edge, this is typically underlined in the context. A connected component is a subset $S$ of vertices such that any two vertices in $S$ are connected, and no vertex in $S^c$ is connected to any vertex in $S$.
With this convention, any two different connected components merge at the minimal connection time of a pair of vertices $(k,l)$ (where $k$ is from one, and $l$ from the other component) to form a single connected component.
Let the mass of any connected component be equal to the number of its vertices.
Due to elementary properties of independent exponentials, it is immediate that a pair of connected components merges at the rate equal to the product of their masses.
In other words, the continuous-time random graph evolves according to the {\em \MC} dynamics:
\begin{equation}
\label{merge}
\begin{array}{c}
\mbox{ any pair of components with masses (sizes) $x$ and $y$ merges }\\
\mbox{ at rate $xy$ into a single component of mass $x+y$.}
\end{array}
\end{equation}
Due to the relation $p= 1-e^{-t}$, this continuous-time random graph exhibits the same phase transition as $G(n,\cdot)$ at time $\approx 1/n$, as $n$ diverges. 

Aldous \cite{aldRGMC} extended this construction as follows: instead of mass $1$, let vertex $i\in [1,n]$ have initial mass $x_i>0$. 
For each $i,j\in [1,n]$ let the edge between $i$ and $j$ appear at rate $x_i x_j$, independently of others. 
The same elementary property of exponentials implies that the transition mechanism is again (\ref{merge}).
Furthermore \cite{aldRGMC}, Proposition 4 shows that if the set of vertices is ${\mathbb N}$, and if $\bx=(x_1,x_2,\ldots) \in l^2$, where the initial mass of $i$ is $x_i$, then this (infinite) random graph process is still well-defined, and its connected component masses form an $l^2$ vector a.s.~at any later time. 
Here is a more precise formulation: let $(\cvd,d)$ be the metric space of infinite sequences
${\bf y} = (y_1,y_2,\ldots)$
with $y_1 \geq y_2 \geq \ldots \geq 0$
and $\sum_i y_i^2 < \infty$,
where 
$d({\bf y},{\bf z}) = \sqrt{\sum_i (y_i-z_i)^2}$.
Let ``${\rm ord}$'' be the ``decreasing
ordering" map defined on infinite-length vectors.
Let $X_i(t)$ be the $i$th largest connected component mass in the above defined random graph.   
The process $(\bX(t),\,t\geq 0)\equiv ((X_1(t),X_2(t),\ldots),t\geq 0)$ started from $\bX(0) = {\rm ord}(\bx) \in \cvd$ is $\cvd$-valued Feller process evolving according to (\ref{merge}) (see \cite{aldRGMC} Propositions 4 and  5, and Section 2.1 in \cite{VLthesis} for an alternative derivation of the Feller property).
Starting with \cite{aldRGMC}, any such process $\bX$ is referred to as a {\em \MC}.
In this note, a {\em graph representation of \MC} (or an {\em MC graph representation} for short) will be any random graph-valued process such that its corresponding ordered component sizes evolve as a \MC.

The processes from the previous two paragraphs are clearly MC graph representations (provided one naturally extends $[1,n]$ to ${\mathbb N}$, and lets $x_i=0$ for $i\geq n+1$).
A different but similar MC graph representation was explored in Bhamidi et al.~\cite{bhamidietal2} (see Section 2.3.1 therein) and recalled next.
Here for each $i,j\in [1,n]$ a (new) directed edge $i\rightarrow j$ appears at rate $x_i x_j/2$, and for each $i$ a self-loop $i \rightarrow i$ appears at rate $x_i^2/2$. The random-graph is effectively an (oriented) multi-graph (a graph with loops and multiple edges allowed).
If the connected components are obtained by taking into account all the edges (regardless of their orientation), and the mass of each connected component is again the sum of masses of its participating vertices,
it is easy to see that the resulting ordered component masses evolve again according to the {\MC} transitions. Indeed, the (potential) presence of multi-edges and loops does not change the connectivity properties or the component masses, so the random graph process from this paragraph can be matched to that from the last paragraph.

As already hinted above,  in this setting it is convenient to embed finite vectors into an infinite-dimensional space.
Refer henceforth to $\bx=(x_1,x_2,x_3\ldots)\in \cvd$ as {\em finite}, if  for some $i\in \mathbb{N}$ we have $x_i=0$. Let the {\em length} of $\bx$ be the number $\len(\bx)$ of non-zero coordinates of $\bx$.
Fix a finite initial configuration $\bx\in \cvd$.
For each $i\leq \len(\bx)$, let $\xi_i$ have exponential (rate $x_i$)
distribution, independently over $i$. 

The order statistics of $(\xi_i)_{i\leq \len(\bx)}$ are denoted by $(\xi_{(i)})_{i\leq \len(\bx)}$. 
Given $\xi$s, define simultaneously for all $q>0$
\begin{equation}
\label{defZbxq}
Z^{\bx,q}(s):= \sum_{i=1}^{\len(\bx)} x_i 1_{(\xi_i/q\, \leq \ s)} -s = \sum_{i=1}^{\len(\bx)} x_{(i)} 1_{(\xi_{(i)}/q\, \leq \ s)} -s, \ s\geq 0, \ q>0.
\end{equation}
In words, $Z^{\bx,q}$ has a unit {negative drift and successive positive jumps, which occur precisely at times $(\xi_{(i)}/q)_{i\leq \len(\bx)}$, and where the $i$th  successive jump is of magnitude $x_{(i)}$.
One could define $Z^{\bx,0}(s):=-s$, but this is not the most natural limit of $Z^{\bx,q}$ as $q\searrow 0$. The point is that, as $q\searrow 0$, $\xi_\cdot/q$ diverge but also that the distance between them diverges, and this latter divergence is more important for the coupling with the random graph.

In \cite{multcoalnew} the family (\ref{defZbxq}) was called the {\em simultaneous breadth-first walks}. It was shown in \cite{multcoalnew}, Proposition 7 that, as $q$ increases, the excursion lengths of the reflected $Z^{\bx,\cdot}$ have the law of the  {\MC} started from the configuration $\bx$. This was an essential step in the proof of  \cite{multcoalnew}, Theorem 2, a strengthening of which is stated as part of Conjecture \ref{Con:one} in Section \ref{S:Scaling}.

A related graph representation of {\MC} was implicitly exhibited in \cite{multcoalnew}, and will be recalled in Section \ref{S:Coupl} below.
In this new setting the random graph process 
is (finite) forest-valued, almost surely. 
Section \ref{S:Coupl} will introduce two extensions, connecting MC forest-valued representations and the continuous-time random graph (defined in the opening paragraphs of this article). 

The coupling with surplus constructions given in Sections \ref{S:take1} and \ref{S:take2} are intrinsic (up to randomization) to the simultaneous breadth-first walk.  
To to the best of author's knowledge, they also carry more detailed information than any of the previous surplus edge studies (compare with \cite{aldRGMC,bhamidietal2,bromar15}). 
In particular, provided that all the labels (positions) are kept for the surplus edges, the continuous-time random graph and the 
`` enriched'' (simultaneous) breadth-first walks are equivalent, either in the sense of the marginal (see Lemma \ref{L:strictatq}) or the full distribution (see Lemma \ref{L:strictatqfull}). 
Moreover, the coupling of Section \ref{S:take2} naturally motivates an extension to the multi-graph setting (linked to \cite{bhamidietal2}) in Section \ref{S:Multi-g}.
In Section \ref{S:Scaling} novel scaling limits are anticipated.

For general background on the random graph and the stochastic coalescence the reader is referred to \cite{aldous_survey, bertoin-fragcoal, bol-book, durrett-rgd,pitman-stflour}, and for specific as well as more recent references to \cite{multcoalnew}.
The edges in this paper will often be defined as oriented, however when the global connectedness in a resulting forest or (multi-)graph is studied, these orientations will not be important.
All the line segments drawn by hand in Figures \ref{F:five}--\ref{F:seven} are meant to be perfectly straight.

\section{The essential surplus and the simple random graph}
\label{S:Coupl}
The following figure illustrates the breadth-first exploration of vertices in a finite rooted tree. 
\begin{figure}[h]
\centering
\includegraphics[scale=0.15]{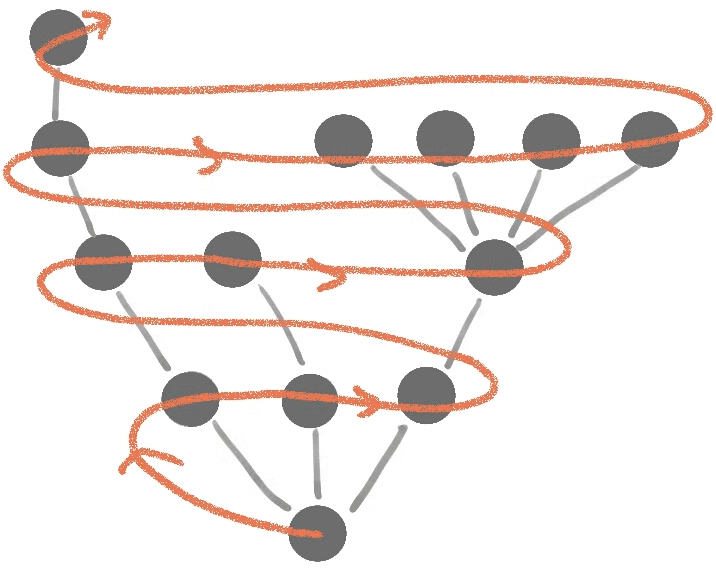}
\caption{The exploration process starts at the root, visits all of its children (these vertices are the 1st generation), then all the children of all the vertices from the 1st generation (these vertices are the 2nd generation), 
and keeps going until all the vertices are visited.}
\label{F:one}
\end{figure}

In this section, random forests will conveniently span the components of a {\em coupled} {\MC}.  
In Section \ref{S:take1} these processes will be explored similarly to \cite{multcoalnew}, with a new feature: the non-spanning or {\em surplus} or {\em excess} edges will be recorded in addition. After that, in Section \ref{S:take2} another graph representation will be proposed in order to preserve the monotonicity.

\subsection{Breadth-first order induced forest and surplus}
\label{S:take1}
Recall (\ref{defZbxq}).
Let us define for $l\leq \len(\bx)$
$$
\pi_l:=i \mbox{ if and only if } \xi_i=\xi_{(l)}, \ l\in [1,\len(\bx)].
$$ 
In this way, $(x_{\pi_1},x_{\pi_2},\ldots, x_{\pi_{\len(\bx)} } )$  is the size-biased random ordering of the initial non-trivial block masses.
Given $a<b$ and an interval $[c,d]$ where $0\leq c <d$, define
$(a,b]\oplus [c,d]:=(a+c,b+d]$, and denote by $|I|$ the Lebesgue measure of a Borel set $I$ (typically $I$ will be an interval below). 
We start by recalling the construction from \cite{multcoalnew} as follows:\\
- set $k=1$ and $F_k^q\equiv F_1^q:=[0,\xi^q_{(1)}]$, $m_k=2$ and let 
$I_1^q=(\xi^q_{(1)}, \xi^q_{(1)} +x_{\pi_1}]$;\\
- for any $l\geq m_k$, 
\begin{equation}
\label{ErecurIla}
\mbox{as long as 
$\xi_{(l)}^q \in   I_{l-1} \mbox{ and } l\leq \len(\bx)$
 define }
I_l^q\equiv I_l:= 
 I_{l-1} 
\oplus [0,x_{\pi_l}], \mbox{  increase $l$ by $1$};
\end{equation}
- let $l_1\equiv l_1(k)$ be equal to thus attained $l$, if  $l_1\leq \len(\bx)$,\\
 increase $k$ by $1$ and let
\begin{equation}
\label{ErecurIlb}
F_k^q:= (\sup(I_{l_1-1})\,,\xi^q_{(l_1)}] \mbox{ and }
I_{l_1}^q\equiv I_{l_1}:= (\xi^q_{(l_1)}, \xi^q_{(l_1)}+x_{\pi_{l_1}}], \ m_k:=l_1+1;
\end{equation}
go back to (\ref{ErecurIla}); \\
- else if  $l_1> \len(\bx)$ (leave $k$ unchanged and) exit.

Note that most of the above quantities depend on $q$, even though the notation does not always indicate it.
The above recursion algorithm naturally produces for each $q>0$ a random (labeled, weighted) forest $\forn(q)$ with vertex set $[1,\len(\bx)]$. 
Extend this definition to a trivial forest of $\len(\bx)$ rooted trees (for example, ordered according to $\pi$) at $q=0$.

For $l=1$ or $l$ equal to some $l_1(q)$ in (\ref{ErecurIlb}), the corresponding $\pi_l$ is the root of a (new) tree in $\forn(q)$, and $I_l^q \setminus  I_{l-1}^q = I_l^q$. For any other $l\leq \len(\bx)$, the corresponding vertex $\pi_l$ is not a root.
For each $l\geq 1$, the dynamics ``listens'' for the children of $\pi_l$ during  $I_l^q\setminus  I_{l-1}^q$.
If $h$ is such that $\xi_{(h)}^q \in I_l^q\setminus  I_{l-1}^q$, then $\pi_l$ (resp.~$\pi_h$) is the parent (resp.~child) of $\pi_h$ (resp.~$\pi_l$) in $\forn(q)$. In symbols $\pi_h \rightarrow \pi_l$ (thinking that any directed edge points to the parent, and therefore towards the root). 

Define 
$$ B^{\bx,q}(s):= Z^{\bx,q}(s) - \inf_{u\leq s} Z^{\bx,q}(u), \ s\geq 0, \ q >0.$$
Due to the above made observations, if (and only if) at time $s$ a new vertex is seen (or `` heard'') in $\forn(q)$, $Z^{\bx,q}$  (and therefore $B^{\bx,q}$) makes a jump up by the amount equal to the mass of that vertex. Upon the completion of the $k$th recursion in (\ref{ErecurIla}), the total sum of upward jumps of $B^{\bx,q}$ during 
$T_k(q)\equiv T_k:=Cl(I_{l_1(k)-1}^q)$
 is entirely compensated by the unit downward drift of $Z^{\bx,q}$ (or $B^{\bx,q}$) during $T_k$. It is also easy to see that $B^{\bx,q}(s)> 0$ in the interior of $T_k$. 
Therefore the $k$th excursion of $(B^{\bx,q},\ q>0)$ above $0$ has the length precisely equal to the sum of the masses of all the blocks (vertices of the tree) explored during $T_k$.
It was initially noted in \cite{multcoalnew} (this approach is rooted in \cite{aldRGMC,EBMC}) that for each fixed $q$, the ordered excursion lengths $(T_l(q),\,l=1,\ldots,k(q))$ have the {\MC} distribution, evaluated at time $q$.
In particular, the exiting $k=k(q)\leq \len(\bx)$ in the above algorithm, clearly equals to the number of connected components at time $q$ of the coupled (continuous time) random graph. 

Proposition 7 in \cite{multcoalnew} shows the following improvement of the just noted distributional equivalence: the excursion lengths $(T_l(q),\,l=1,\ldots,k(q))$ have the random graph (or \MC) law as a process in $q$, where the initial state at time $0$ is $\bx$.
One might wish to strengthen this in saying that the 
 above coupling of $Z^{\bx,q}$ and $\forn(q)$ provides a bijective matching between the $k$th excursion of $B^{\bx,q}$ above $0$ (necessarily started at  $\xi_h^q$ for some $h\in [1,\len(\bx)]$), and a spanning (breadth-first search) tree (rooted at $\pi_h$) of the unique component of the continuous-time random graph which contains $\pi_h$  at time $q$. 

Some care is however needed here.
Indeed, while the mergers of different connected components (or different subtrees of $\forn(\cdot)$) arrive at precisely the multiplicative (or random graph) rate, the new edges arriving in $\forn$ that correspond to those mergers always connect the root of one of the components (its excursion starts later in $B^{\bx,\cdot\, -}$) to the last visited (or listed) vertex in the other component (its excursion is the one starting just before in $B^{\bx,\cdot \, -}$).
In addition, within each connected component  the edges evolve according to an unusual ``prune and reconnect'' rule, where vertices and subtrees are gradually ``moved closer'' to the root. 
In particular, the forest $(\forn(q),\,q\geq 0)$ is not a monotone process with respect to the order induced by the subgraph relation.
Still, $\forn(\cdot)$ is a random graph (forest) valued process, whose tree masses evolve precisely according to (\ref{merge}), so it is an MC graph representation. 
We refer to any edge of $\forn(q)$ as {\em spanning edge at time} $q$. 

Furthermore,  there is a natural way to build the surplus or excess edges on top of $\forn(q)$, for each $q$ separately, so that the distribution of the resulting graph is exactly that of the continuous-time random graph at time $q$.
We need some additional notation in order to explain this extension.
For each $h,k \in [1,\len(\bx)]$ such that $k\geq h$ , let $\zeta^{h,k}$ be a Poisson process of marks arriving at rate $x_k$, independently over $h$ and $k$. 
Since there are no loops or multi-edges in the continous-time random graph $G(\len(\bx),1-e^{-q}),\,q\geq 0)$, we only need to check if there is a surplus edge between each pair of vertices which
are connected, but not by an edge in $\forn(q)$.

Suppose that for some $h\in [2,\len(\bx)]$ we have $\xi_{\pi_h}^q\in I_{h-1}^q$.
If $l>h$ is such that $\xi_{\pi_l}^q\in I_{h-1}^q$ (or equivalently, if $\pi_l$ is another vertex, attached to the same tree after $\pi_h$ but prior to any child of $\pi_h$) then (due to the breadth-first exploration order) there are only two different possibilities:\\
(a) $\pi_l$ is of the same generation as $\pi_h$, or\\
(b) $\pi_l$ is from the next generation (then necessarily a child of some vertex $\pi_k$ with $k<h$ in the same generation as $\pi_h$).\\
Indeed, all the other vertices $\pi_k$ (with $k>h$) that belong to the same tree in $\forn(q)$ could not have been heard before the start of $I_h^q \setminus I_{h-1}^q$. See Figure \ref{F:two} for an illustration. 
\begin{figure}[h]
\centering
\includegraphics[scale=0.15]{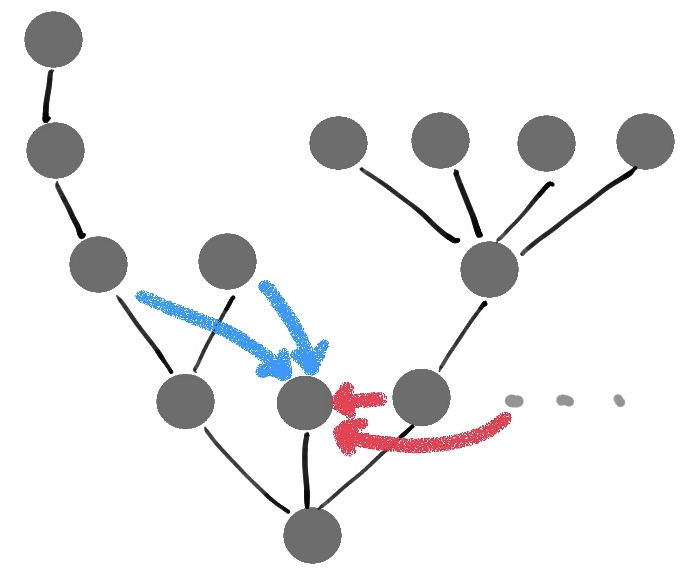}
$\hspace{1cm} $
\includegraphics[scale=0.15]{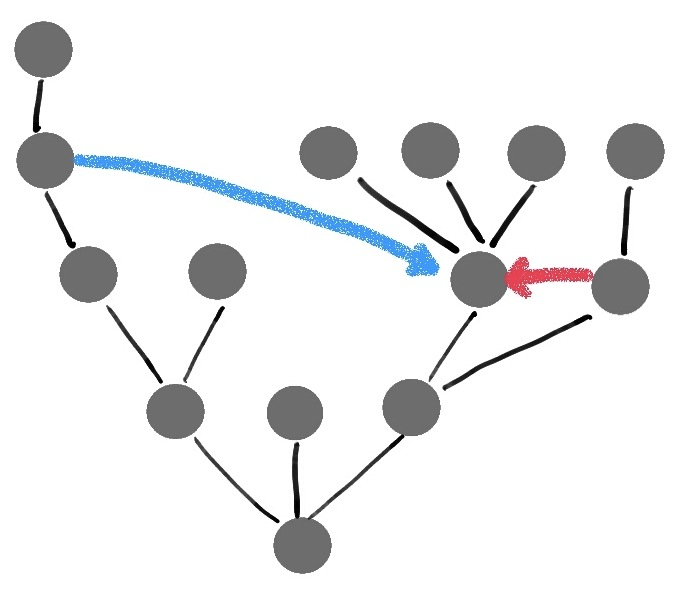}
\caption{Two different configurations are shown with $h=3$ and $h=7$ respectively. Possible surplus edges coming from vertices in case (a) are indicated in red, and those from vertices in case (b) are indicated in blue.}
\label{F:two}
\end{figure}

Figure \ref{F:three} below shows the surplus ``influence region''  for a typical non-root vertex.
\begin{figure}[h]
\centering
\includegraphics[scale=0.12]{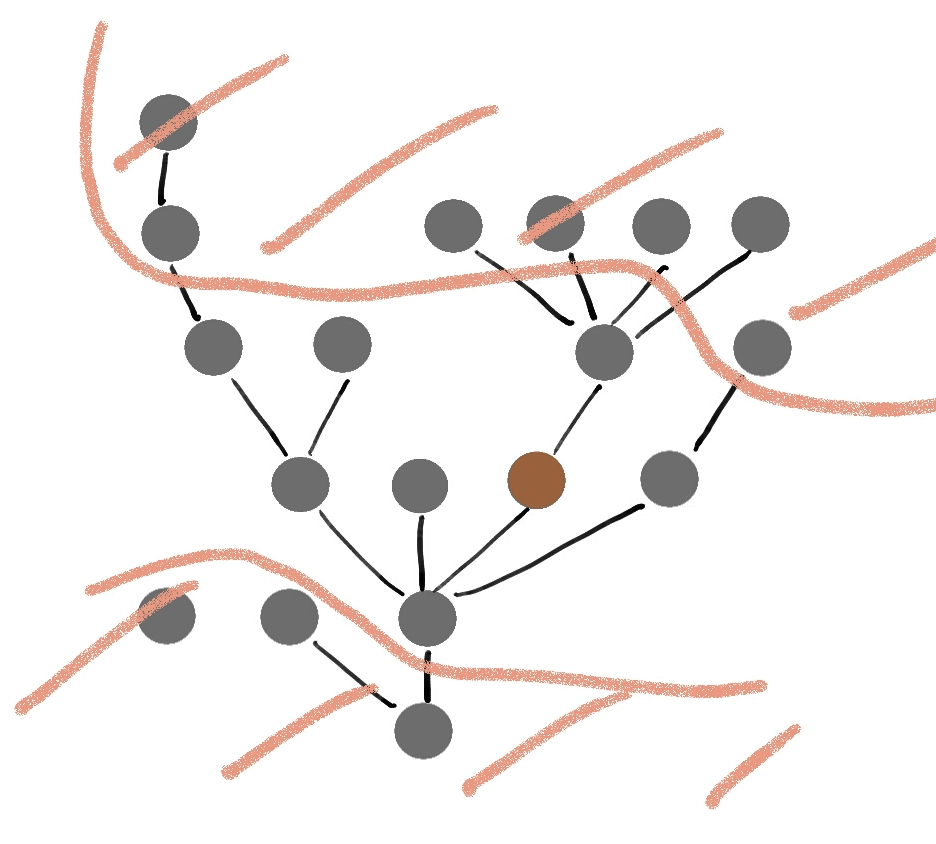}
\caption{Due to the breadth first order, the brown vertex can be connected by a surplus edge only to a vertex outside of the barred regions.}
\label{F:three}
\end{figure}

During $I_h^q\setminus I_{h-1}^q$ consider those (and only those) $\zeta^{h,l}$ with $l>h$ and $\xi_{\pi_l}^q\in I_{h-1}^q$. For any such $l$ draw the (red or blue) surplus edge connecting $l$ and $h$ if and only if  $\zeta^{h,l}[0, q x_{\pi_h}] \geq 1$ (or equivalently  iff $\frac{1}{q}\zeta^{h,l}[0,x_{\pi_h}]\geq 1$). It is clear that this edge appears with probability $1-e^{-q x_{\pi_l} x_{\pi_h}}$, independently of everything else. Each spanning edge at time $q$ is also present with the probability $1-e^{-q m_1 m_2}$, where $m_1$ and $m_2$ are the masses of the two vertices (see above, or Section 2 in \cite{multcoalnew} for more details). 
Denote by $G(q)=(V,E(q))$ the resulting random graph, where $E(q)$ is the union of the spanning and the surplus edges at time $q$.
One can record the just made observations as follows.
\begin{Lemma}
\label{L:strictatq}
For each $q \geq 0$, $G(q)$ has the law of the continuous-time random graph evaluated at time $q$.
\end{Lemma}


Suppose that one is only interested in counting the surplus edges in various connected components of the random graph (without keeping track of their exact position).
Then note that
the joint (total) intensity of all the Poisson marking processes $\zeta^{h,l}$ active during $I_h^q \setminus I_{h-1}^q$ is 
\begin{equation}
\label{Etotalinten}
q \cdot \sum_{l>h,\, \xi_{\pi_l}^q\in I_{h-1}^q} x_{\pi_l} \equiv q\left(\sum_{s\in  I_{h-1}^q} \Delta Z^{\bx, q}(s) - |I_{h}^q| \right).
\end{equation}
The quantity on the right hand side is precisely $q (B^{\bx,q}(\sup(I_{h-1}^q))- x_{\pi_h})$, a rescaled value (up to an error) of the reflected $Z^{\bx,q}$ at the right endpoint of $I_{h-1}^q$ (recall that $Z^{\bx,q}$ attains its infimum on $[0,\sup(I_h^q)]$ at the left endpoint $\inf(I_h^q)$ of $I_h^q$).

\subsection{Monotone forest representation and surplus}
\label{S:take2}
We next revise and enrich (via an additional randomization) the coupling algorithm just described. 
Instead of the random forest process $(\forn(q),\, q>0)$ another forest-valued process $(\foro(q),\, q>0)$ is described, so that:\\
(i) for each $q$ the (sub)trees in $\forn(q)$ and in $\foro(q)$ are ordered in the same way, and consist of exactly the same (random) subsets of vertices, almost surely, however\\
(ii) $\forn(q)$ and $\foro(q)$ typically (and for large $q$ very likely) have different edge sets,
(iii) $\foro(q_1)$ is a subgraph of $\foro(q_2)$ whenever $q_1\leq q_2$, almost surely.\\
Furthermore, in this new setting the surplus edges will again be accounted for in a rather natural but somewhat more complicated way (in comparison to that of Section \ref{S:take1}). 
 
The initial forest $\foro(0)$ is again trivial, and therefore equal to $\forn(0)$.
During a strictly positive (random) interval of time $\foro(\cdot)$ will remain $\foro(0)$.
At some (stopping) time $Q_1>0$ the first connection is established between $\pi_{L_1}$ and $\pi_{L_1-1}$ such that $\xi_{(L_1)}^{Q_1}\in I_{L_1-1}^{Q_1}$ (and necessarily $\xi_{(h)}^{Q_1} \not \in I_{h-1}^{Q_1}$ for any $h\neq L_1$, 
$h\in \{2,\ldots, \len(\bx)\}$).
At time $Q_1$ both $\forn$ and $\foro$ make the same jump: the new edge $\pi_{L_1}\rightarrow \pi_{L_1-1}$ appears. 
After that, during an interval of time of positive (random) length, $\foro$ stays equal to $\foro(Q_1)$, and eventually a new connection occurs at some time $Q_2>Q_1$. 
The difference between $\foro$ and $\forn$ may be visible already at time $Q_2$. 
Indeed, suppose that it is $\xi_{(L_1+1)}^{Q_2}$ that enters $I_{L_1}^{Q_2}$ at this very moment. As already noted, in  $\forn(Q_2)$ we must have $\pi_{L_1+1}  \rightarrow \pi_{L_1}$ at time $Q_2$. 
However, in $\foro(Q_2)$ the new edge is either\\
- $\pi_{L_1+1}  \rightarrow \pi_{L_1}$, with probability $x_{\pi_{L_1}}/(x_{\pi_{L_1-1}}+ x_{\pi_{L_1}})$, or\\
- $\pi_{L_1+1}  \rightarrow \pi_{L_1-1}$, with the remaining probability $x_{\pi_{L_1-1}}/(x_{\pi_{L_1-1}} + x_{\pi_{L_1}})$.\\
In $\foro$ this (and any other) parent and child relation, once established, will stay fixed throughout the evolution. As already noted, in $\forn$ both $\pi_{L_1}$ and $\pi_{L_1+1}$ will eventually become children of the same vertex $\pi_{h}$ for some $h\in \{1,\ldots, L_1-1\}$ and this $h$ is bound to change until finally becoming equal to $1$.

The complete construction of $\foro$ is as follows: $\foro(0)=\forn(0)$; for $q>0$, $\foro(q)=\foro(q-)$ unless $q$ is such that the number of components (trees) in $\forn$ at time $q$ decreases by $1$.
The latter happens if and only if for some $i \in [1,\len(\bx)]$, $\xi_{\pi_i}^{q-} \not \in I_{i-1}^{q-}$ ($\pi_i$ is a root in $\forn(q-)$, or equivalently, in $\foro(q-)$) and $\xi_{\pi_i}^{q} \in I_{i-1}^{q}$. For such $q$, let $\foro(q)$ inherit all the edges of $\foro(q-)$, and in addition draw a new edge 
\begin{equation}
\label{Enewedge}
\pi_i \rightarrow \pi_{L(q)} \mbox{ in } \foro(q),
\end{equation}
where $L(q)$ is chosen at random (in the size-biased way) from the vertices of the ``tree to the left of $\pi_i$'' in $\foro(q-)$. More precisely, conditionally on $\{\foro(s),\, s<q\}$ and $\{\forn(s),\, s\leq q\}$, apply (\ref{Enewedge}) with $L(q)$ equal to $l$ with probability $x_{\pi_l}/|I_{i-1}^{q-}|$, for each $l$ such that $\xi_{\pi_l}^{q-} \in I_{i-1}^{q-}$ (or equivalently, for each $l$ such that $\pi_l$ and $\pi_{i-1}$ are connected in $\foro$ and in $\forn$ during $(q-\eps,q)$ for all sufficiently small positive $\eps$).

\smallskip
It is clear from the just presented construction that $\foro$ is a monotone forest-valued process, and also that, almost surely, for each $q$ the trees in $\foro(q)$ are composed of exactly the same vertices as the (corresponding) trees in $\forn(q)$ (or equivalently, $i$ and $j$ are connected in $\foro(q)$ if and only if they are connected in $\forn(q)$).  In short, the above properties (i)--(iii) apply to $\foro$, hence $\foro$ is another MC graph representation, realized on the same probability space as $\forn$ and $Z^{\bx,q}$.

{\em Remark.}
The just described way of attaching edges in $\foro$ is the most natural (if not the only) choice for a monotone MC forest representation coupled with the breadth-first walk. Indeed, in order to respect the order of vertices (and connected components) which is induced by the coupled walks $Z^{\bx,\cdot}$, one has no option but to attach the root of the tree to the right (in $\foro(q-)$) to some vertex of the tree to the left, and (\ref{Enewedge}) means that the parent $L$ is picked uniformly from the mass in the component (tree) to the left.

\subsubsection{Surplus on top of $\foro$}
\label{S:surtopforo}
We refer to any edge of $\foro(q)$ as {\em spanning edge arriving before} $q$. Note the difference with the similar (and weaker) definition in terms of $\forn$. 

It is next described how to account for the surplus edges in a way compatible to the coupled breadth-first walks $Z^{\bx,\cdot}$ (this was also a feature of the construction 
in Section \ref{S:take1}). The present goal is to preserve the monotonicity of the surplus edges in time. As the reader will see in the next section, there is a natural multi-graph that emerges from this construction, and it happens to be the random multi-graph from Bhamidi et al.~\cite{bhamidietal2}, recalled in the introduction.

Consider the following figure, showing a (part of a) realization of $Z^{\bx,q}$, with the ``space under the curve'' of the corresponding $B^{\bx,q}$  split into conveniently chosen polygons (to be soon split further into parallelograms), and the triangles marked by different colors for increased readability.
\begin{figure}[h]
\centering
\includegraphics[scale=0.25]{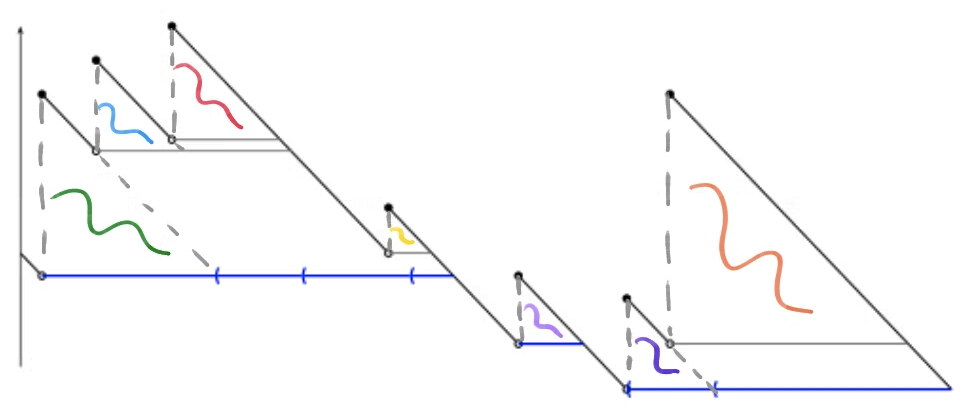}
\caption{The middle excursion is still in its original (and simplest possible) state - it corresponds to the tree consisting of a single vertex. The excursion to the right corresponds to a tree with two vertices, and the one to the left corresponds to a tree with four vertices. The dashed lines are added to indicate the triangles, the $i$th of which is to be matched to vertex $\pi_i$.}
\label{F:four}
\end{figure}

For $q$ close to $0$ the curve $Z^{\bx,\cdot}$ has only $\len(\bx)$ many (moving) {\em triangular} excursions. As $q$ increases, merging gradually happens, and simultaneously (due to the coupling described in previous sections) the excursions of $Z^{\bx,\cdot}$ get more complex.
It is interesting to describe here the exact structure of this {\em excursion mosaic} process, induced by the gradual ``pile up'' of the original $\len(\bx)$ triangles on ``top of each other''. 

The excursion mosaic (or mosaic) is made by drawing at each time $q$ the horizontal (blue in Figure 4) segment on the basis of each excursion of $Z^{\bx,q}$. Furthermore, if the merging of a pair of trees with roots $\pi_i$ and $\pi_j$ ($i<j$) occurs at time $q$, then\\
- the blue segment starting at $\xi_i^q$ is extended to the end of the new excursion at time $q$;\\
- the blue segment starting at $\xi_j^q$ turns gray, and is lifted at any later time $q+z$ to the vertical level $Z^{\bx,q+z}(\xi_j^{q+z}-)$;\\
- for each of the original $\len(\bx)$ triangles, the grey line indicating its hypotenuse is extended below, until it meets the  unique horizontal blue segment of the corresponding excursion.

From now on, call any excursion of $Z^{\bx,q}$, with its corresponding blue and grey segments obtained according to the above given  procedure, {\em ornamented}.
For an ornamented excursion of $Z^{\bx,q}$ which contains $\xi_{(k)}^q$, we say that {\em it carries} $\pi_k$.
It should be clear that the just given mosaic construction  and the related definitions can be transposed so that, for each $q$, the reflected processes $B^{\bx,q}$ has the same ornamented excursions  as $Z^{\bx,q}$ (with the difference that in $B^{\bx,q}$ the excursions start from the abscissa).

If at time $q>0$ there is a non-trivial spanning tree in $\forn$ (or equivalently in $\foro$) of length four, its corresponding ornamented excursion must have one of the five forms depicted in Figure \ref{F:five}.
\begin{figure}[h]
\centering
\includegraphics[scale=0.08]{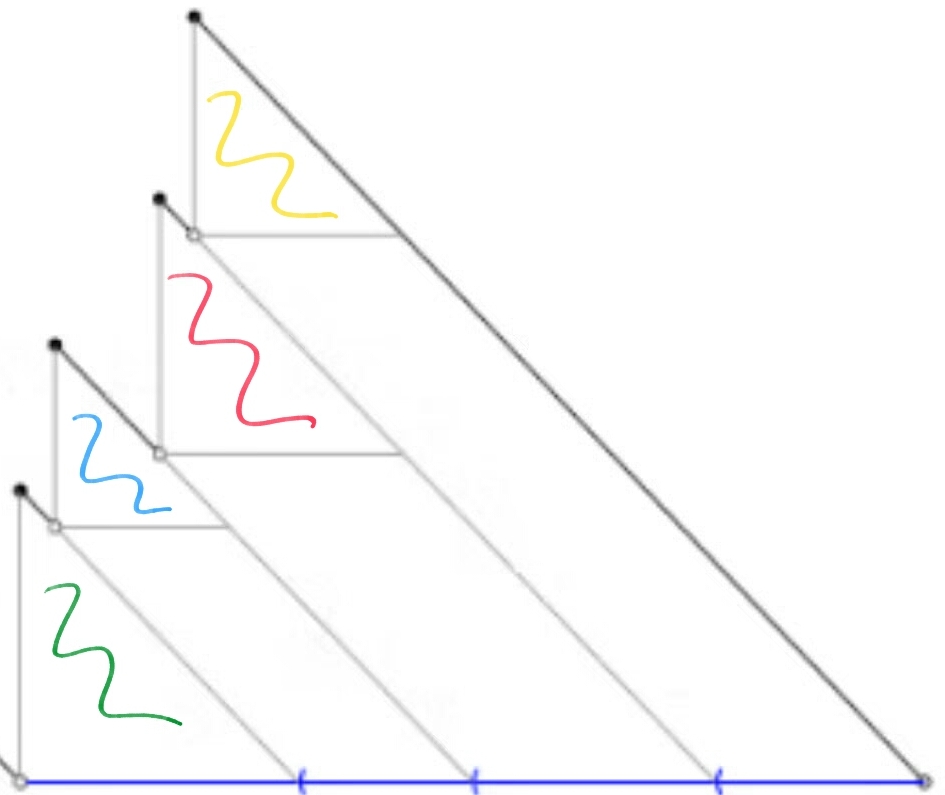} \
\includegraphics[scale=0.08]{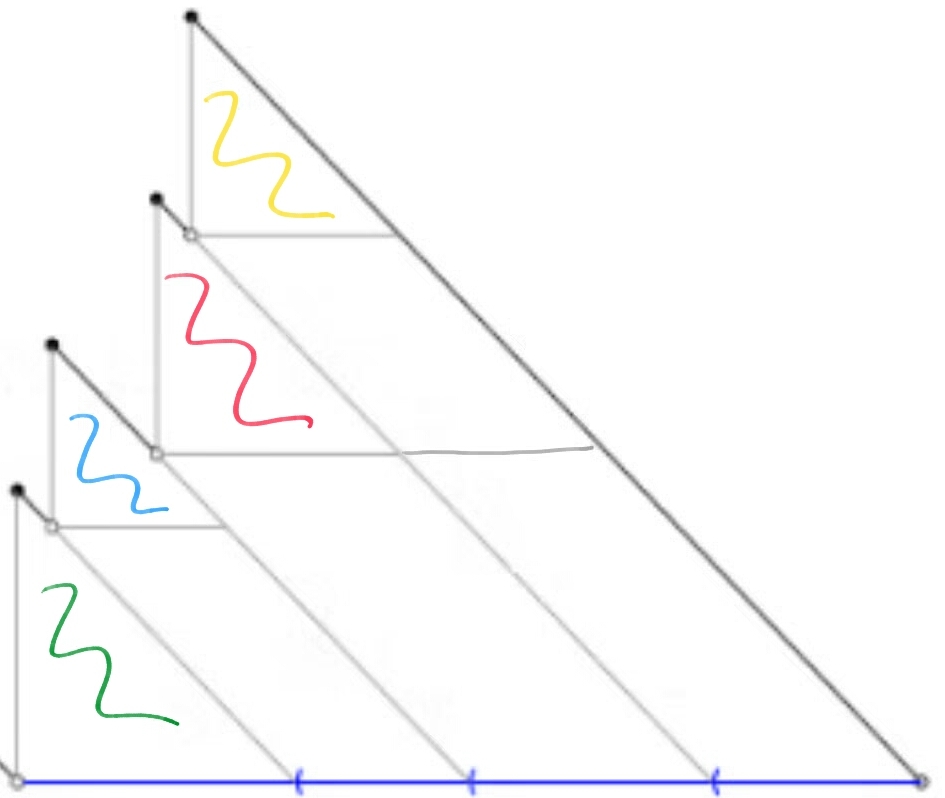} \
\includegraphics[scale=0.08]{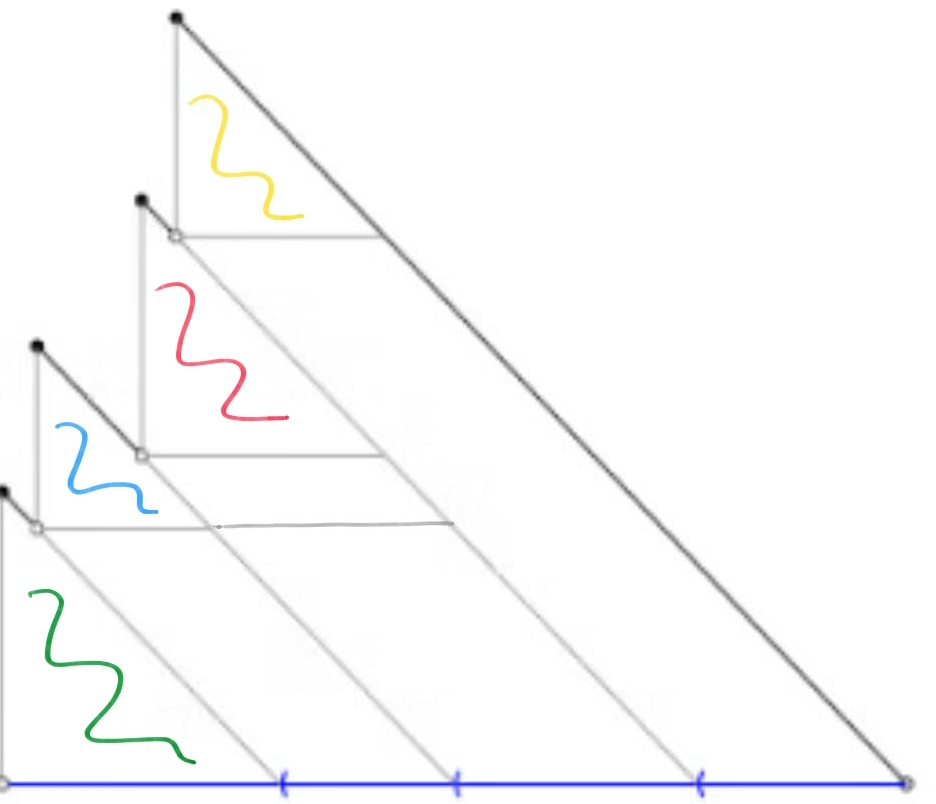} \
\includegraphics[scale=0.08]{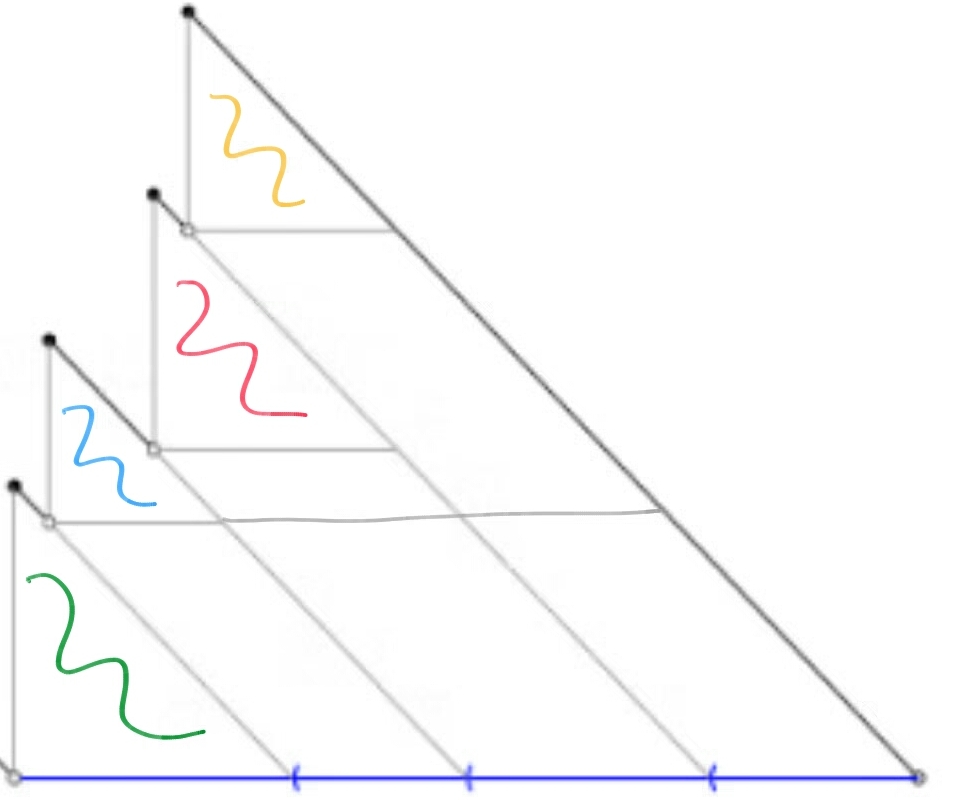} \
\includegraphics[scale=0.08]{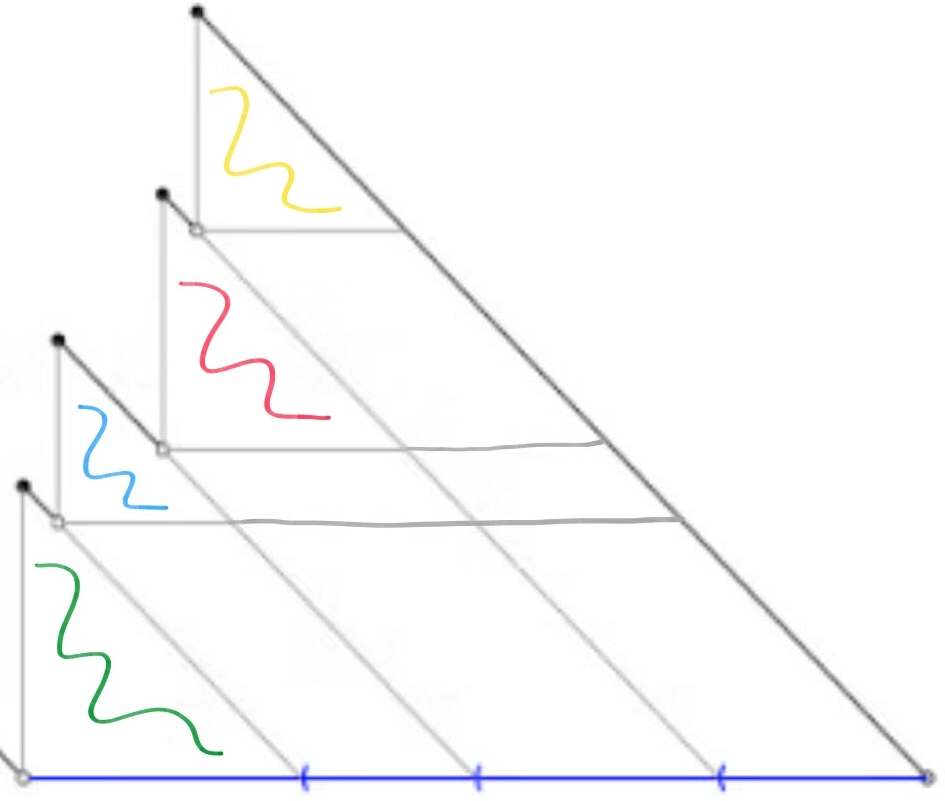}
\caption{Possible realizations of the excursion mosaic, restricted to a single excursion, viewed at a fixed time $q>0$. The tree corresponding to the fifth image is completely specified by the mosaic, while for the other figures there is at least one edge picked at random according to (\ref{Enewedge}). For example, in the third image the yellow vertex may be connected to either the green, blue, or red one, and in the fourth image it can be connected to either the blue or the red vertex.}
\label{F:five}
\end{figure}

The following figure shows the ornamented excursions from the first two images in Figure \ref{F:five} with their corresponding trees, as well as some ill-defined ornamentations.
\begin{figure}[h]
\centering
\includegraphics[scale=0.08]{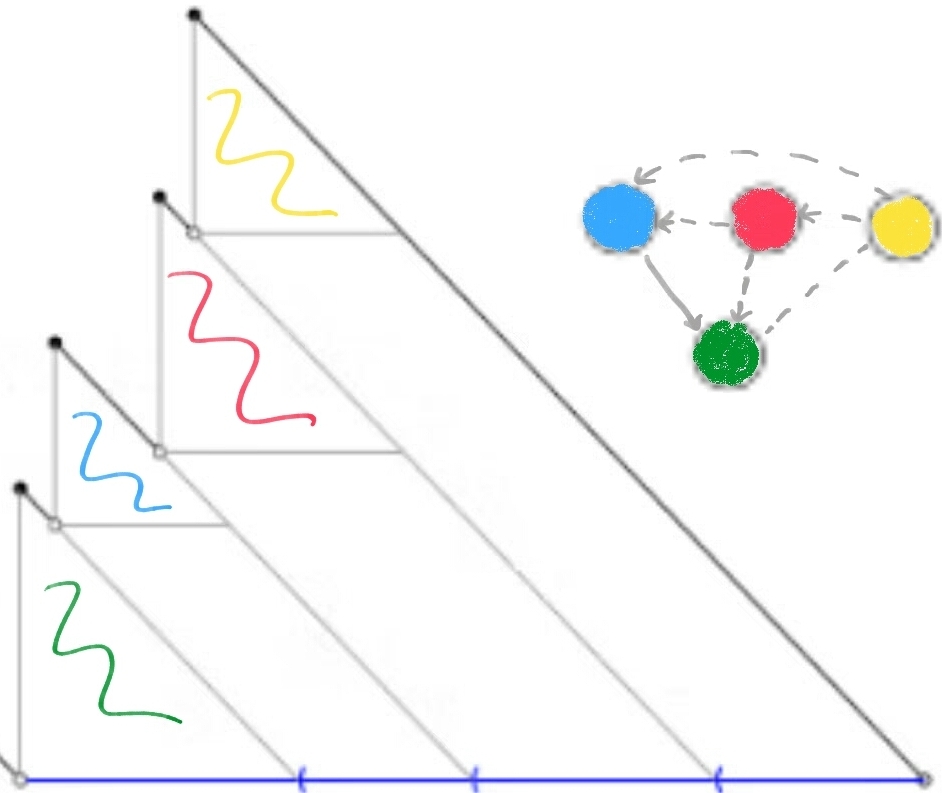} \ \ 
\includegraphics[scale=0.08]{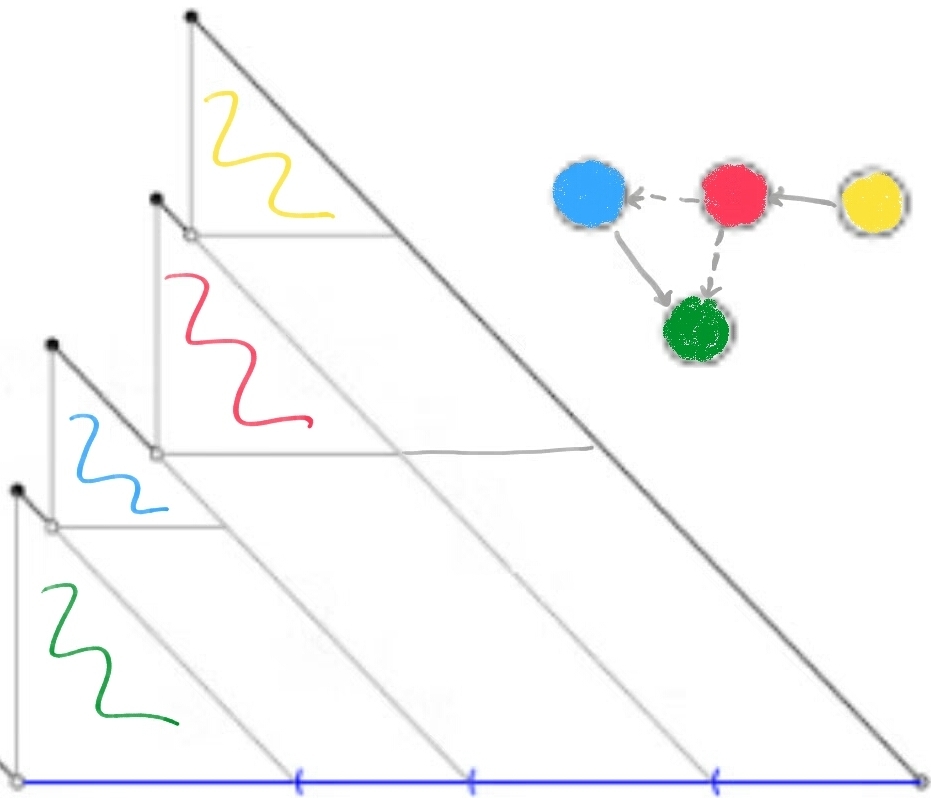} \  \ \
\includegraphics[scale=0.08]{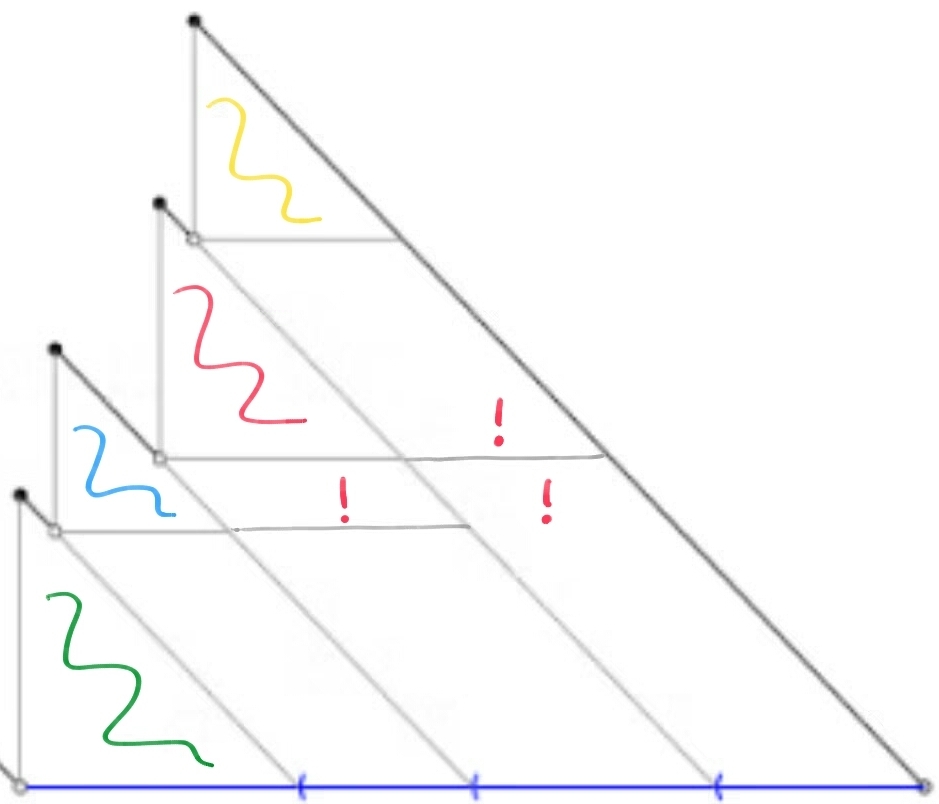} \ \
\includegraphics[scale=0.08]{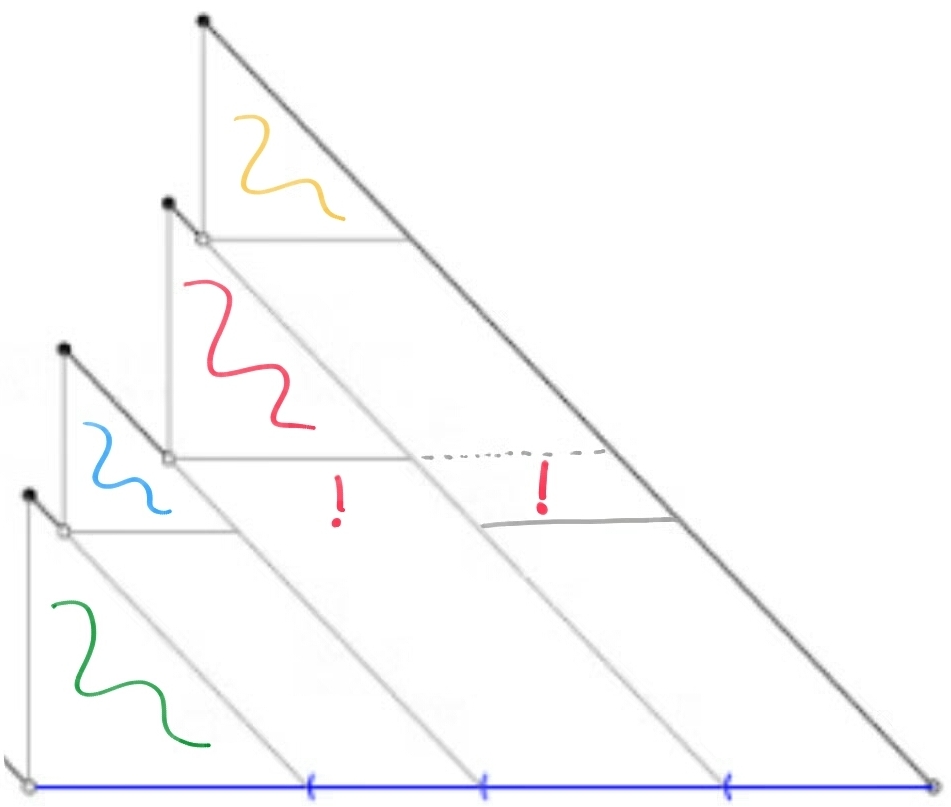}
\caption{A dashed arrow indicates a possible edge (from each non-root vertex there is a unique directed edge issued), while a full arrow indicates an edge specified by the mosaic.
In the final two images the exclamation marks indicate joint presence or absence of grey segments, which are impossible in an excursion mosaic.}
\label{F:six}
\end{figure}

Let us first show  on an example how the surplus edges can be superimposed in a consistent way. Consider the second image in Figure \ref{F:six}. Enumerate the four parallelograms specified by the mosaic in some way, for example traverse them row-by-row from left to right, as an analogue to the breadth-first order (see Figure \ref{F:one}). Let $\zeta^{b;g}$, $\zeta^{r;b,g}$, $\zeta^{y;b,g}$ and $\zeta^{y;r}$ be independent Poisson point processes, matched respectively to these four regions. As $q$ increases, for each of the parallelograms the base stays fixed in length (although it moves closer to the origin), and the height increases. 
Let a point arrive to $\zeta^{b;g}$ at rate $x_b \cdot x_g$, to $\zeta^{r;b,g}$ at rate $x_r (x_b + x_g)$, to $\zeta^{y;b,g}$ at rate $x_y (x_b + x_g)$, and to $\zeta^{y;r}$ at rate $x_y \cdot x_r$.\\
As it turns out, most of these Poisson point processes will not be needed now (they are nevertheless defined with an intention of later use).
Let us assume that, due to randomization (\ref{Enewedge}), the realization of the corresponding random tree has the following three edges:  $b \rightarrow g$, $r \rightarrow g$,  and  $y \rightarrow r$. 
We only need to watch for marks in $\zeta^{r;b,g}$ and $\zeta^{y;b,g}$.
When a new mark arrives to $\zeta^{r;b,g}$, 
with probability $x_b/(x_b+x_g)$ the process creates an edge $r \rightarrow b$, unless this edge already exists. Nothing happens with the remaining probability.
When a new mark arrives to $\zeta^{y;b,g}$,
with probability $x_b/(x_b+x_g)$ (resp.~$x_g/(x_b+x_g)$) the process creates an edge $y \rightarrow b$ (resp.~$y \rightarrow g$),  unless it already exists.
It is likely clear to the reader that these transitions are chosen with the purpose of preserving the random graph transitions within the connected components.

\smallskip
In the general case, 
 one has the collection of Poisson processes $\zeta^{l;j-k}$, where $j,k,l$ run over all the indices in $[1,\len(\bx)]$ such that
\begin{equation}
\label{Econstr} 
j\leq k\leq l, \mbox{ and if }k=l\mbox{ then also }j=l.
\end{equation} 
The processes $(\zeta^{l;l-l})_l$ play a special role, to be explained in Section \ref{S:Multi-g}.  
If $l>k$, then $\zeta^{l;j-k}$ is in charge of generating a surplus edge from $\pi_l$ to some vertex in the range $\{\pi_j,\ldots,\pi_k\}$, but {\bf only when compatible} with the excursion mosaic. More precisely, $\zeta^{l;j-k}$ will be {\em active} starting from  time $T_{l;j-k}$ at which the ornamented excursion of $Z^{\bx,\cdot}$ carrying $\pi_l$ 
 merges with an ornamented excursion carrying precisely $\pi_j,\ldots,\pi_k$. 
In this way, the random time $T_{l;j-k}$ depends on the mosaic. 
It can happen that $j=k<l$, if and only if the ornamented excursion carrying $\pi_l$ merges with a simple triangular excursion carrying $\pi_j$ only.
On the event $\{T_{l;j-k} = \infty\}$, the corresponding $\zeta^{l;j-k}$ is never activated.
On $\{T_{l;j-k} < \infty\}$, the behavior of $\zeta^{l;j-k}$ after $T_{l;j-k}$ is a generalization of the one given in the example above (see e.g.~$\zeta^{r;b,g}$).
More precisely, the points arrive to $\zeta^{l;j-k}$ at rate $x_{\pi_l} (x_{\pi_j} + \ldots + x_{\pi_k})$. 
The surplus edges are created as follows: given an arrival to $\zeta^{l;j-k}$ at time $q$, an index $I$ is drawn (independently from the past of the mosaic, of $\foro$, and of the surplus edge data), so that $I=i$ with probability $x_{\pi_i}/( x_{\pi_j} + \ldots + x_{\pi_k})$ for each $i \in \{j,\ldots,k\}$. Given $I$, the surplus edge $l\rightarrow I$ is created, unless it already exists. 
For each $q>0$, call any edge created in this way before time $q$ a {\em surplus edge arriving before} $q$.
Let $G_1(q)=(V,E_1(q))$, where $E_1(q)$ is the union of the spanning and surplus edges arriving before $q$.
Then it is clear that $G_1$ is a monotone random graph process: $G_1(q_1) \subset G_1(q_2)$ whenever $q_1\leq q_2$. In fact we have a stronger claim, that may serve as a main motivation for the extra construction presented in this section.
\begin{Lemma}
\label{L:strictatqfull}
Let $\len(\bx)=n$. If $x_1=x_2=\ldots=x_n=1$, the law of $(G_1(q),\,q\geq 0)$ equals that of $(G(n,1-e^{-q}), q\geq 0)$.
For general initial weights $x_1,\ldots, x_n$,   
$(G_1(q),\,q\geq 0)$
is a realization of the Aldous' (inhomogenous) continuous-time random graph, as recalled at the beginning of the second introductory paragraph. 
\end{Lemma}

\section{Unlimited surplus and a canonical multi-graph}
\label{S:Multi-g}
Here we focus on the second construction above (see Section \ref{S:surtopforo}).
In particular, recall the excursion mosaic, and the family of compatible Poisson point processes
$\zeta^{l;j-k}$, where $j,k,l$ satisfy the constraints given in (\ref{Econstr}).
For each $l$, $\zeta^{l;l-l}$ should now be matched at time $q$ to the triangle spanned by the points 
$(\xi_{(l)}^q, Z^{\bx,q}(\xi_{(l)}^q-))$, $(\xi_{(l)}^q, Z^{\bx,q}(\xi_{(l)}^q))$
and
$(\xi_{(l)}^q+x_{\pi_l}, Z^{\bx,q}(\xi_{(l)}^q-))$ (or equivalently, to that spanned by $(\xi_{(l)}^q, B^{\bx,q}(\xi_{(l)}^q-))$,  $(\xi_{(l)}^q, B^{\bx,q}(\xi_{(l)}^q))$ and $(\xi_{(l)}^q+ x_{\pi_l}, B^{\bx,q}(\xi_{(l)}^q-))$). This Poisson point process is active already at time $0$ (the triangular excursions exist from the very beginning). 
We therefore define $T_{l;l-l}\equiv 0$ almost surely, for each $l\in [1,\len(\bx)]$. 
Points arrive to $\zeta^{l;l-l}$ at rate $x_{\pi_i}^2/2$, and at the time of each arrival, a self-loop $\pi_l \rightarrow \pi_l$ is created.

{\em Remark.}
The factor of $1/2$ is natural if one thinks of each original block as continuous ``spread'' of mass, and of each self-loop as a connection between two points on the block. Provided that these sites are sampled independently and uniformly at random (as is always done in the Erd\"os-Ren\`yi setting), 
there are two possible ways of choosing the same loop. 
Not surprisingly, this factor is also natural from the perspective of matching the total surplus edge count to the area under the curve (excursion) of $Z^{\bx,\cdot}$ (or of $B^{\bx,\cdot}$), to be explained soon.
\hfill\endrem

As in Section \ref{S:surtopforo}, for $l>k$ the counting process $\zeta^{l;j-k}$ is activated at time $T_{l;j-k}$ (hence never on  $\{T_{l;j-k}=\infty\}$). After activation, the points arrive to $\zeta^{l;j-k}$ again at rate  $x_{\pi_l} (x_{\pi_j} + \ldots + x_{\pi_k})$.
The surplus multi-edges are created as before, but without an additional ``lack of previous presence'' constraint: given an arrival to $\zeta^{l;j-k}$ at time $q$, an index $I$ is drawn in the same way as before, and a new surplus (multi-)edge $l\rightarrow I$ is created at time $q$. The thus obtained multi-graph is a version of that from \cite{bhamidietal2}.

It was already explained how $\zeta^{l;l-l}$ can be matched to the $l$th triangular region under the curve $B^{\bx,\cdot}$.
It is useful to make explicit here that $\zeta^{l;j-k}$ (on $\{T_{l;j-k} < \infty\}$) can analogously be matched to a parallelogram shaped region (evolving in time) on the mosaic, for any consistent choice of $l>k\geq j$. 
Before time $T_{l;j-k}$ this parallelogram does not exist, exactly at time $T_{l;j-k}$ it has height $0$, and its height (strictly) increases at any future time. 
Indeed, this parallelogram of constant base length $x_{\pi_l}$ is created at time $T_{l,j-k}$ by the excursion mosaic, and at any time $z>T_{l,j-k}$ it is specified by the four lines
\begin{equation}
\label{Elines}
\begin{array}{lc}
y= - (x - \xi_{(l)}^z) +Z^{\bx,z}(\xi_{(l)}^{z}-) & \ \mbox{\small{"left" boundary}},\\
y= Z^{\bx,z}(\xi_{(j)}^{z}-) & \ \mbox{\small{"bottom" boundary}},\\
y= - (x - \xi_{(l)}^z) +Z^{\bx,z}(\xi_{(l)}^{z}) & \ \mbox{\small{"right" boundary}},\\
y= Z^{\bx,z}(\xi_{(j_1)}^{z}-) & \ \mbox{\small{"top" boundary}},\\
\end{array}
\end{equation}
where $j_1$ is the minimal index larger than $j$ (hence necessarily inside $[k+1,l]$) such that $T_{l,j_1-k_1}<T_{l,j-k}$, for some $k_1\in [j_1,l]$. 

{\em Remark}.
A careful reader will note (or easily derive from the definitions) that $j_1$ is in fact (almost surely on $\{T_{l;j-k}<\infty\}$) equal to $k+1$. \hfill \endrem

Figure \ref{F:seven} shows how the mosaic drawn in previous figures might look at a later time. 
For the sake of illustration, let us assume that the eight jumps (vertices) in the figure 
are $\pi_3,\ldots,\pi_{10}$, where $\pi_3$ is indicated in green, and $\pi_{10}$ in blue.

\begin{figure}[h]
\centering
\includegraphics[scale=0.15]{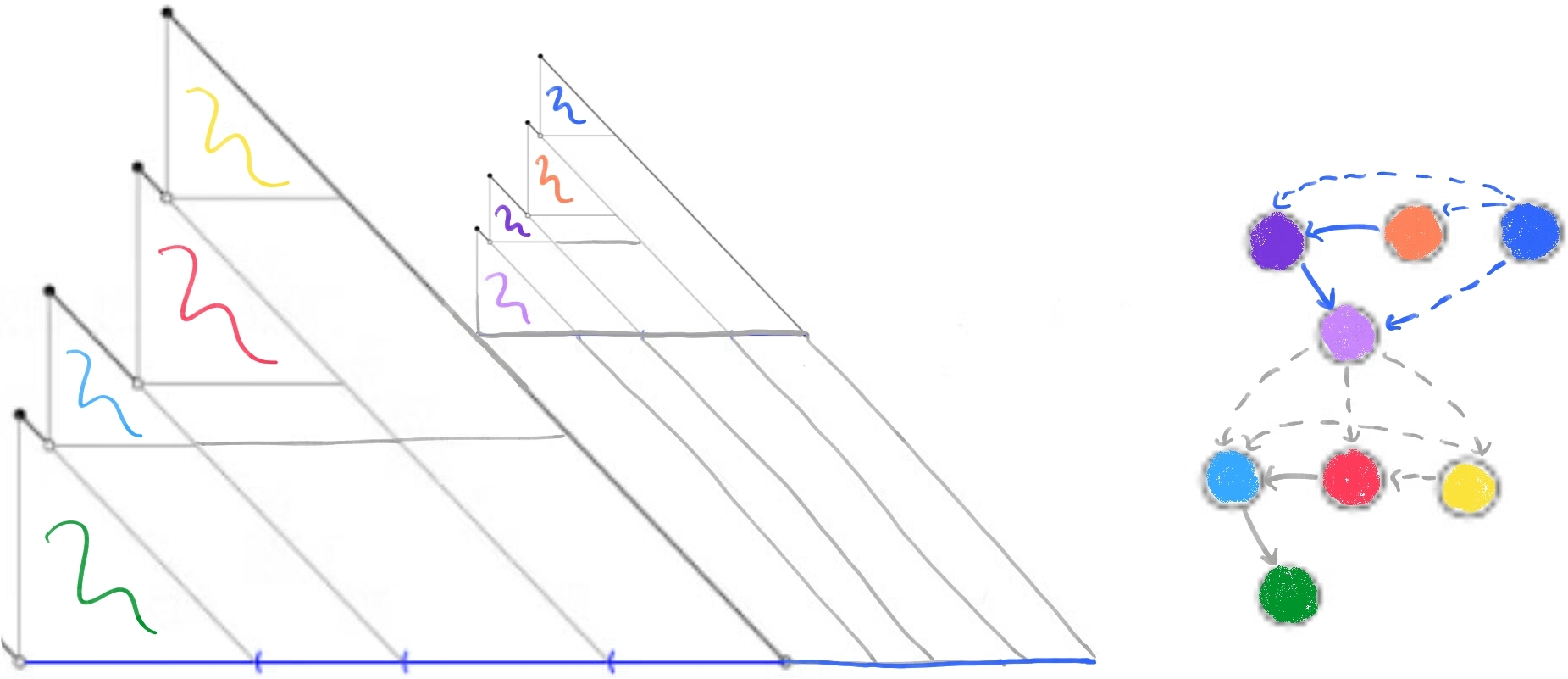} 
\caption{The tree is again partially determined by the mosaic, the dashed edges indicate various options, to be decided according to (\ref{Enewedge}).
For the (orange) vertex $\pi_9$, there are three active Poisson point processes $\zeta^{9;8-8}$, $\zeta^{9;7-7}$ and $\zeta^{9;3-6}$ running at that time.
For the (blue)  $\pi_{10}$, there are only two active Poisson point processes $\zeta^{10;7-9}$, and  $\zeta^{10;3-6}$ running. 
The process $\zeta^{9;7-8}$ corresponds to the middle parallelogram in the ``$\pi_9$-slice'', containing the triangle marked in orange, and $\zeta^{10;3-6}$ corresponds to the lower parallelogram in the ``$\pi_{10}$-slice'', containing the triangle marked in blue. 
There is no $\zeta^{7;3-5}$, why? Find $\zeta^{6;3-5}$ and $\zeta^{6;4-5}$.
}
\label{F:seven}
\end{figure}
The excursion mosaic is an object of interest due to the following claim in particular.
\begin{Proposition}
\label{P:rates}
For any choice of $j,k,l$ specified in (\ref{Econstr}), and for any fixed $q$, on $\{T_{l;j-k}\leq q\}$,
the cumulative arrival rate to $\zeta^{l;j-k}$ up to time $q$ equals (almost surely) the area of the region 
 in the excursion mosaic matched to $\zeta^{l;j-k}$ at time $q$, multiplied by $q$.
\end{Proposition}
{\em Proof.} If $j=k=l$, the statement is trivial to check.
Suppose that $l>k$, hence the corresponding region in the mosaic is a parallelogram. 
Then clearly the total arrival rate to $\zeta^{l;j-k}$ before time $q$ equals $(q-T_{l;j-k}) \cdot x_{\pi_l} (x_{\pi_j} + \ldots + x_{\pi_k}).$.
Abbreviate $L:=x_{\pi_j} + \ldots + x_{\pi_k}$.
Then the total rate above can be rewritten as 
$$
q \left(L- \frac{L T_{l;j-k}}{q}\right) \cdot x_{\pi_l}.
$$ 
The claim is therefore that the height of the parallelogram specified in (\ref{Elines}) equals the term in parentheses.

Recall the remark following (\ref{Elines}). 
At time $q$, the height of the parallelogram from (\ref{Elines}) is $Z^{\bx,q}(\xi_{(k+1)}^{q}-)-Z^{\bx,q}(\xi_{(j)}^{q}-)$.
It was noted from the start that $T_{l;j-k}$ (whenever finite) is the time of merging of two excursions (one carrying $\pi_l$ and another carrying exactly $\pi_j,\ldots,\pi_k$ just before that time). 
On $T_{l;j-k}<\infty$ the excursion carrying $\pi_l$ prior to time $T_{l;j-k}$ must have $\pi_{k+1}$ as its initial vertex, or equivalently, the root of the corresponding tree in \foro\ must be $\pi_{k+1}$ (see also the most recently mentioned remark). 
So the merging of these two excursions can be written out as identity
$$
\xi_{(k+1)} - \xi_{(j)} = L T_{l;j-k}.
$$
It is simple to check now that $Z^{\bx,q}(\xi_{(k+1)}^{q}-)-Z^{\bx,q}(\xi_{(j)}^{q}-)$ is almost surely equal to 
$$
L- \xi_{(k+1)}^q - \xi_{(j)}^q \equiv L- \frac{\xi_{(k+1)} - \xi_{(j)}}{q},
$$
so the identity required above is true.
\hfill \endpf

\medskip
We shall refer to the bounded (random) region specified by the lines 
\begin{equation}
\label{Elines_one}
\begin{array}{lc}
y= - (x - \xi_{(l)}^z) +Z^{\bx,z}(\xi_{(l)}^{z}-) & \ \mbox{\small{"left lower" boundary}},\\
x= \xi_{(l)}^{z}, & \ \mbox{\small{"left upper" boundary}}\\
y= - (x - \xi_{(l)}^z) +Z^{\bx,z}(\xi_{(l)}^{z}) & \ \mbox{\small{"right" boundary}},\\
y= 0 & \ \mbox{\small{"bottom" boundary}},\\
\end{array}
\end{equation}
as the {\em $\pi_l$-slice (under the curve of $B^{\bx,q}$)}.
From the discussion above we see that the $\pi_l$-slice is split by the mosaic into the right isosceles triangle of area $(x_{\pi_l})^2/2$ and (possibly) additional parallelograms, each of which corresponds to $\zeta^{l;j-k}$ for some $j\leq k<l$ such that $T_{l;j-k}\leq q$.
In this way, the region under the excursion of $B^{\bx,z}$ carrying exactly $\pi_h,\pi_{h+1},\ldots, \pi_{h+r}$ at time $q$ is the union of $\pi_h$-slice, $\pi_{h+1}$-slice, $\ldots$, and $\pi_{h+r}$-slice.
The intersection of adjacent regions and slices are sets (segments) of zero Lebesgue measure.
We therefore obtain
\begin{Corollary}
\label{C:rates}
Almost surely for each $q>0$, \\
(a) the cumulative rate of (oriented) surplus edges issued from $\pi_l$ before time $q$ is the area of $\pi_l$-slice at time $q$, multiplied by $q$,\\
(b) the cumulative rate of surplus edges in the component of $G_1(q)$ consisting of $\pi_h,\pi_{h+1},$ $\ldots, \pi_{h+r}$ is the area of the excursion of $Z^{\bx,q}$ (or of $B^{\bx,q}$) carrying $\pi_h,\pi_{h+1},\ldots, \pi_{h+r}$, multiplied by $q$.
\end{Corollary}

\section{Scaling limits: a discussion}
\label{S:Scaling}
Both expression (\ref{Etotalinten}) and Corollary \ref{C:rates} are promising in view of novel scaling limits for near-critical random graphs, outside of the domain of attraction of the Aldous standard \MC.
At present the only (eternal) augmented \MC\ is the original one of Bhamidi et al.~\cite{bhamidietal2}, a version of which was constructed in \cite{bromar15}  as the scaling limit of the random graph with surplus counts for special initial configurations of the form $x_1=x_2=\ldots= x_n=1/n^{2/3}$, $0=x_{n+1}=\ldots$, as $n$ diverges. 

Let $\kappa\geq 0$, $\tau\in \mathbb{R}$, and $\bc \in l^3$ with non-increasing components. 
Given a family $(\xi_j')_j$ of independent exponentials, where $\xi_j'$ has rate $c_j$,
define 
\begin{equation}
 V^\bc (s) = \sum_j \left(c_j 1_{(\xi_j' \leq s)} - c_j^2s \right)
, \ s \geq  0 . \label{defVc}.
\end{equation}
For each $t\in \mathbb{R}$, let
$$
W^{\kappa,t- \tau, \bc}(s) = \kappa^{1/2}W(s) -\tau s - \sfrac{1}{2}\kappa s^2  + V^\bc (s) + t s, \ s \geq 0 \label{defWtc},
$$
where
$W$ is standard Brownian motion, and $W$ and $V^\bc$ are independent, and let
\begin{equation}
 B^{\kappa, t-\tau,\bc}(s) := W^{\kappa, t-\tau,\bc}(s) - \min_{0 \leq s^\prime \leq s} W^{\kappa, t-\tau,\bc}(s^\prime), \ s \geq 0. \label{defBtc}
\end{equation}
For each $t$, let $\bX(t)=\bX^{\kappa,\tau,\bc}(t)$ be the infinite vector of ordered excursion lengths of $B^{\kappa, t-\tau, \bc}$ away from $0$. 
Theorem 2 in \cite{multcoalnew} says that (for most parameters $\kappa,\tau,\bc$) $(\bX(t), t\in (-\infty,\infty))$ is a realization of the extreme eternal
\MC\ corresponding to $(\kappa, \tau, \bc)$. 
Let $N$ be a homogeneous Poisson point process on $[0,\infty)\times [0,\infty)$, independent of $\sigma\{W,V^\bc\}$.
In analogy to \cite{aldRGMC,bhamidietal2}, let $N^{\kappa, t-\tau,\bc}(s)$ be the number of points in $N$
 below the curve $u\mapsto B^{\kappa,t-\tau,\bc}(u)$, $u\in [0,s]$. 
To each excursion of $B^{\kappa, t-\tau, \bc}$ above $0$, one can assign a random ``mark count'' to be the increase in $N^{\kappa, t-\tau,\bc}$  attained during this excursion (see \cite{bhamidietal2}, Section 2.3.2 for details). 
Let $Y_i(t)$ be this count corresponding to the $i$th longest excursion of $B^{\kappa,t-\tau,\bc}$, and $\bY(t)=(Y_1(t),Y_2(t),\ldots)$.

Given the observations made in previous sections, the following can be anticipated: 
\begin{Conjecture}
\label{Con:one}
For each parameter triple $(\kappa,\tau,\bc)$, as in \cite{EBMC,multcoalnew}, Theorem 2,  $({\bZ^{\kappa,\tau,\bc}}(t)=(\bX(t),\bY(t))$, $ \,-\infty < t< \infty)$ is an  {\em eternal  augmented \MC} corresponding to $(\kappa,\tau,\bc)$. Furthermore, ${\bZ^{\kappa,\tau,\bc}}$ is the simultaneous scaling limit of near-critical random graph component sizes and surplus counts,  under the hypotheses of the initial configurations given in  \cite{EBMC}, Proposition 7. The extreme eternal augmented \MC s are only the constant ones, and the non-trivial ones given here (corresponding to valid parameters $(\kappa,\tau,\bc)$). Any eternal augmented \MC\ is a mixture of extreme ones. 
\end{Conjecture}

The excursion mosaic and the accompanying PPP family $\zeta^{l;j-k}$, $j\leq k\leq l$ (see Section \ref{S:surtopforo}) has a much richer structure than the mere component sizes superimposed by surplus edge counts. 
Is there a natural framework and candidate for its scaling limit in the near-critical regime(s)? 
This insight would surely encompass a clearer understanding of mark counts $\bY$  in the eternal augmented coalescents $\bZ$.

\end{document}